\documentclass{amsart}
\usepackage[T1]{fontenc}
\usepackage{amssymb}
\usepackage[hyperref]{degt}

\def\2{\color{red}}

\def\minitab#1{\vcenter\bgroup\rm
 \def\minitabskip{#1}%
 \def\-##1{\setbox0\hbox{$00$}\hbox to\wd0{\hss$##1$\hss}}%
 \def\.{\cdot}
 \let\\\cr
 \halign\bgroup\strut\hss$##$&&#1\hss$##$\hss\cr}
\def\endminitab{\crcr\egroup\egroup}

\def\trule{\noalign{\vspace{1pt}\hrule\vspace{2pt}}}
\def\HHH{\setbox0\hbox{$(0)$}\hbox to\wd0{\hss$\HH$\hss}}
\def\counts#1{\multispan2\strut\hss#1:\ignorespaces}

\def\clist#1{\omit\minitabskip$#1$\hss}
\def\tabconfig#1{\hbox to12pt{\hss$#1$\!\!\hss}}

\def\dt{\mathinner{\ldotp\ldotp}}

\def\HH{\Sigma}
\def\cross{\underline}
\def\cross#1{\setbox0\hbox{$#1$}\mathord{\vrule height2.8pt depth-2.4pt width\wd0\relax
 \kern-\wd0\relax\box0}}

\def\maxall{M}
\def\numall{\#}
\def\maxC{\maxall_\C}
\def\maxR{\maxall_\R}
\def\maxRR{\maxall_\R^\circ}
\def\numC{\numall_\C}
\def\numR{\numall_\R}
\def\numRR{\numall_\R^\circ}
\let\subtext\mathrm
\let\subtext\text
\def\spec{_\subtext{spec}}
\def\nspec{_\subtext{ns}}
\def\both{_\subtext{both}}
\def\hyp{_\subtext{hyp}}
\def\units{^{\times\!}}
\def\girth{\operatorname{girth}}
\def\Fn{\operatorname{Fn}}

\def\CN{\Cal N}
\def\CK{\Cal K}
\def\bL{\bold L}
\def\Fano{\Cal{F}}
\def\graph{\Gamma}
\def\genus{\operatorname{genus}}

\title[Real split hyperplane sections]
{Real split hyperplane sections\\ on smooth polarized $K3$-surfaces}

\author{Alex Degtyarev}

\address{%
Department of Mathematics\\
Bilkent University\\
06800 Ankara, TURKEY}

\email{
degt@fen.bilkent.edu.tr}

\thanks{%
The author was partially supported by the T\"{U}B\DOTaccent{I}TAK grant 123F111
\unskip
}



\begin{document}

\begin{abstract}
We find upper bounds, sharp in most cases, on the number of real hyperplane
sections of real smooth polarized $K3$-surfaces that split into lines.
Most bounds coincide with their complex counterparts.
\end{abstract}

\maketitle

\section{Introduction}\label{S.intro}

In recent paper~\cite{degt:hplanes} we obtain sharp upper bounds
$\maxC(2d)$ on the the
number of \emph{split}, \ie, completely split into a union of lines,
hyperplane sections of a smooth polarized $K3$-surface $X\into\Cp{d+1}$ of
degree $h^2=2d$. For short, we call these split hyperplane sections
\emph{$h$-fragments} in the \emph{graph of lines}, or \emph{Fano graph}
$\Fn X$ of~$X$, \latin{a.k.a}.\ the dual adjacency graph of all lines on~$X$. In fact,
for all degrees $2d\ge6$, we obtain the complete list of
\emph{$h$-configurations}, \ie, unions of $h$-fragments in $\Fn X$.
(In particular, $h$-fragments do not exist if $2d>28$.)
This
would allow us to answer many further questions about the geometry of
$h$-fragments.

One of such natural questions is the maximal number of \emph{real}
$h$-fragments on a \emph{real} smooth polarized $K3$-surface. (Recall that a
\emph{real} algebraic variety is a complex one equipped with a
\emph{real structure}, \ie, anti-holomorphic involution $\iota\:X\to X$.)
We need to distinguish between two notions:
\roster*
\item
\emph{totally real} $h$-fragments, \ie, those composed of real lines, and
\item
real $h$-fragments with pairs of $\iota$-conjugate lines allowed.
\endroster
The maximal numbers of, respectively, totally real, real, and complex
$h$-fragments on surfaces of degree~$2d$ are denoted by
\[*
\maxRR(2d)\le\maxR(2d)\le\maxC(2d).
\]
The principal result of the present paper is the following theorem, where we
also cite the complex bounds $\maxC$ from~\cite{degt:hplanes};
furthermore,
to emphasize the similarity, the very few real bounds that differ
from their complex counterparts are underlined.

\theorem[see \autoref{proof.main}]\label{th.main}
The bounds $\maxRR(2d)\le\maxR(2d)\le\maxC(2d)$ are as follows\rom:
\[*
\let\u\underline
\minitab\quad
h^2=2d          &     2&  4& 6& 8&10&12&14& 16&18&20&22&24&28\\\trule
\maxC           &    72& 72&76&80&16&90&12& 24& 3& 4& 1& 1& 1\\
\maxR           &\u{66}&\u?&76&80&16&90&12&\u6& 3& 4& 1& 1& 1\\
\maxRR          &\u{66}&\u?&76&80&16&90&12&\u4& 3& 4& 1& 1& 1\\
\endminitab
\]
In the unsettled case $2d=4$ of quartics in~$\Cp3$ we have
\[*
\underline{50}\le\maxRR(4)\le\maxR(4)\le\underline{59}.
\]
\endtheorem

We conjecture (\cf. \cite[Conjecture~3.3]{degt:hplanes})
that $\maxRR(4)=\maxR(4)=\underline{50}$.

After a brief arithmetical background in \autoref{S.prelim},
\autoref{th.main} is proved in \autoref{proof.main}, where we consider each
degree $2d$ separately and, in each case, provide a number of extra details
such as the values taken by the numbers
\[*
\numRR(X)\le\numR(X)\le\numC(X)
\]
of $h$-fragments on an individual surface~$X$, counts and bounds itemized by
the types of the graphs, \etc.
(Recall that one of the discoveries of~\cite{degt:hplanes} was the fact that,
in most degrees, even for $2d=6$, the combinatorial structure of the graph of
lines constituting an $h$-fragment is not unique.)

One special case that is worth mentioning here is that of octics, $2d=8$.
A very general octic $K3$-surface $X\into\Cp5$ is a
\emph{triquadric}, \ie, regular complete intersection of three quadric
hypersurfaces. However, the space of smooth octics has a divisor of
\emph{special} ones, \viz. those not representable in that form. According
to~\cite{degt:lines}, the configurations of lines on special octics differ
dramatically from those on triquadrics, and so do the $h$-configurations.
Thus,
the bounds $\maxRR(8)=\maxR(8)=\maxC(8)=80$ in \autoref{th.main} are realized
by triquadrics; for special octics we have (see \autoref{s.h=8})
\[
\maxRR(8\spec)=\underline{42},\quad
\maxR(8\spec)=\underline{54}\quad
\text{(\vs. $\maxC(8\spec)=72$, see \cite{degt:hplanes})}.
\label{eq.8spec}
\]

Another interesting case is that of sextics. As recently discovered in
\cite{degt.Rams:sextics}, here we also have a meaningful dichotomy between
\emph{special} and \emph{non-special} sextics, the former being those lying
in a \emph{singular} quadric in $\Cp4$. Similar to octics, special sextics
are detected homologically by the presence of a $3$-isotropic vector in
$\NS(X)$, and this property affects the configurations of lines,
see \cite[Corollary 6.9]{degt.Rams:sextics}.
This time, the bounds $\maxRR(6)=\maxR(6)=\maxC(6)=76$
in \autoref{th.main} are realized by the
\emph{special} configuration $\Psi_{42}$, whereas for non-special sextics
we have
\[
\maxRR(6\nspec)=\maxR(6\nspec)=\maxC(6\nspec)=36,\quad
 \text{realized by $\Theta_{36}$ in~\cite{degt.Rams:sextics}}.
\label{eq.6nspec}
\]
\emph{Unlike} octics, the $h$-fragment $K(3,3)$ may appear in both special
and non-special sextics, so that one may ask about the
configurations realized by both. We have
\[
\maxRR(6\both)=\maxR(6\both)=\maxC(6\both)=16,\quad
 \text{Humbert sextic in~\cite{DDK}}.
\label{eq.6both}
\]

In \autoref{S.hyper} we discuss yet another special case, \viz.
hyperelliptic projective
models.

\subsection{Notation}\label{s.notation}
We use $\bA_n$, $\bD_n$, $\bE_n$ for the \emph{negative definite} root
lattices generated by the simply laced Dynkin diagrams of the same name.

Unless stated otherwise, all lattices in this paper are even and all
algebraic varieties are over the field~$\C$ of complex numbers.

The inline notation $[a,b,c]$ stands for the
rank~$2$ lattice $\Z u+\Z v$, $u^2=a$,
$v^2=c$, $u\cdot v=b$. In particular, $\bU:=[0,1,0]$ is the
\emph{hyperbolic plane}.

We use $L(n)$ for the lattice~$L$ with the form scaled by $n\in\Q$; this is
opposed to the shortcut $nL:=L^{\oplus n}$. By a slight abuse, $-L:=L(-1)$.

All other notations are either standard or explained in the text. Worth
mentioning is the more functorial notation $\discr L:=L^\vee\!/L$ for the
\emph{discriminant group}; it is denoted by $q_L$ in
the founding paper~\cite{Nikulin:forms}.

\subsection{Acknowledgements}\label{s.thanks}
I am grateful to the organizers of the
``ddg40 : \emph{Structures alg\'{e}briques et ordonn\'{e}es}'' conference
(August 4--8, 2025, Banyuls-sur-mer) for the opportunity to present
the results of~\cite{degt:hplanes}. I extend my gratitude to the audience of
the conference for their patience and a number of interesting questions,
including those that inspired the present paper.

\section{Preliminaries}\label{S.prelim}
We adopt the traditional paradigm and mostly treat a $K3$-surface
$X\to\Cp{d+1}$ as its $2d$-polarized N\'eron--Severi lattice $\NS(X)\ni h$,
$h^2=2d$. Furthermore, we assume that $\NS(X)$ is generated by~$h$ and the
classes of the lines on~$X$, at least over $\Q$. For a concise description of
the algorithms involved, we refer to \cite{degt:hplanes} (very brief) or
\cite{degt.Rams:sextics} (extended). For the present paper we merely recall
that, in order to single out an equilinear stratum with a given
configuration~$\graph$ of lines, regarded as a (multi-)graph, we need to
consider tha lattice
\[*
\Fano(\graph):=(\Z\graph+\Z h)/\ker,\quad h^2=2d,\quad
 v^2=-2,\ v\cdot h=1\ \text{for $v\in\graph$},
\]
and fix
\roster*
\item
a finite index extension $N\supset\Fano(\graph)$ admitting a primitive
isometry
\[*
N\into\bL:=2\bE_8\oplus3\bU\simeq H_2(X;\Z),
\]
\item
a representative $T\in\genus(N^{\perp\!})$,
known as the \emph{transcendental lattice} of~$X$,
of the genus
of $N^\perp\subset\bL$
(by~\cite{Nikulin:forms}, $\genus(N^{\perp\!})$ is determined by~$N$), and
\item
a
bijective anti-isometry
\[
\Gf\:\discr N\to\discr T.
\label{eq.phi}
\]
\endroster
For each configuration~$\graph$ appearing in~\cite{degt:hplanes} all these
data are readily available as they are used as part of the classification of
configurations. Therefore, below we take them as an input and concentrate on
detecting the real structures on~$X$ or, more precisely, on a representative
of the corresponding equilinear stratum. Note that we do \emph{not} make an
attempt to discuss the connectedness of the real strata: according
to~\cite{Itenberg,Josi,DI:empty}, this may be much more involved than the
complex case.

\subsection{Real configurations}\label{s.real}
The following statement is
an immediate consequence of the global Torelli theorem for $K3$-surfaces:
given $N,T,\Gf$ as in~\eqref{eq.phi}, the real structures
\emph{on appropriate representatives} of the respective equilinear stratum
are the involutive elements of the set
\[
\OG_{\pm h}^-(N)\times_\Gf\OG^-(T),
\label{eq.real}
\]
where $\OG^-\subset\OG$ is the set of autoisometries reversing the
\emph{positive sign structure}, \ie, coherent orientation of maximal positive
definite subspaces.
Here, the ``appropriate representative'' means that an element
$g\in\eqref{eq.real}$ should be fixed first, whereupon $X$ is chosen
very general and such that its period $\Go\in T\otimes\C$ satisfies
$g^*(\Go)=\bar\Go$.

The finite group $\OG_{\pm h}(N)\supset\OG_{\pm h}^-(N)$ is easily computed:
\[*
\OG_{\pm h}\bigl(\Fano(\graph)\bigr)=\Aut\graph\times\{\pm\id\},\qquad
\OG_{\pm h}(N)=\stab\CK,
\]
where $\CK\subset\Fano(\graph)$ is the kernel of the finite index extension
$N\supset\Fano(\graph)$. For $\Aut\graph$, we use the \texttt{digraph}
package in \GAP~\cite{GAP4.13}. (Needless to say that the bulk of the other
computation is also performed using \GAP~\cite{GAP4.13}.)

The computation of $\OG(T)\supset\OG^-(T)$ is also immediate if $\rank N=20$
(the so-called \emph{singular $K3$-surfaces}), \ie, $T$ is a positive
definite lattice of rank~$2$. In general, finding all involutions in
$\OG^-(T)$ may be more involved, \cf. \cite{Nikulin:sublattice}.
Fortunately, apart from the totally real case in \autoref{s.totally.real}
below,
we need this group only when
\[
\text{$T\simeq2\bU(3)$ is a rescaling of the unimodular lattice $2\bU$}.
\label{eq.2U}
\]
The
involutions
$g\in\OG^-(T)=\OG^-(2\bU)$ are found as in~\cite{Nikulin:forms},
according to their invariant sublattices $g^+:=\Ker(g-\id)$, which may be
\[*
[2],\ \bU,\ \bU(2),\ [2]\oplus[-2],\ \bU\oplus[-2].
\]
(Recall that $g\in\OG^-(T)$ imposes $\Gs_+(g^+)=1$.)
Since the natural homomorphism
\[*
\OG\bigl(2\bU(3)\bigr)\to\Aut\bigl(\discr2\bU(3)\bigr)
\]
is surjective,
for~\eqref{eq.real} we easily find the three (out of
eight) conjugacy classes of involutions
in $\Aut\bigl(\discr2\bU(3)\bigr)$ represented by involutions
$g\in\OG^-\bigl(2\bU(3)\bigr)$.

\subsection{Totally real configurations}\label{s.totally.real}
For the study of totally real $h$-fragments we can confine ourselves to the
case where \emph{all} lines are real: indeed, we can merely replace~$N$ with
its primitive sublattice rationally generated by~$h$ and real lines.

The following criterion (\cf. \cite[Lemma~3.8]{DIS}) is based on
Nikulin's classification of real structures on
$K3$-surfaces~\cite{Nikulin:forms}:
{\em a finite index extension $N\supset\Fano(\graph)$ is realized by a real
polarized $K3$-surface with all lines $\ell\in\graph$ real if and only if
there exists a transcendental lattice $T\in\genus(N^{\perp\!})$ such that}
\[
T\supset[2]\quad\text{or}\quad T\supset\bU(2).
\label{eq.totally.real}
\]
In other words, we need a finite index extension of $N\oplus[2]$ or
$N\oplus\bU(2)$ in which $N$ is primitive and which admits a primitive
isometry to~$\bL$.
Using \cite{Nikulin:forms},
denoting by $r:=\rank T=22-\rank N$ the rank of
$T$ and letting $\CN:=\discr N$ (and
$\CN_p:=\CN\otimes\Z_p$ for a prime~$p$),
this criterion can be restated as follows.

A transcendental lattice $T\supset[2]$ exists if and only if
\roster*
\item
for each prime $p>2$, one has either $\ell(\CN_p)\le r-2$ or $\ell(\CN_p)=r-1$
and, in the latter case,
$\det\CN_p=-2\ls|\CN|\bmod(\Z_p\units)^2$;
\item
either $\ell(\CN_2)\le r-2$ or $\ell(\CN_2)=r$ and, in the letter case,
there is a vector $u\in\CN_2$ of order~$2$ and square $-\frac12\bmod2\Z$
(hence, necessarily an orthogonal direct summand)
such that either
\roster*
\item
$u$ is not characteristic, $u\cdot v\ne v^2\bmod\Z$ for some $v\in\CN_2$ of
order~$2$, or
\item
$u$ is characteristic and $\det u^\perp=\pm2\ls|\CN|\bmod(\Z_2\units)^2$.
\endroster
\endroster
\emph{Assuming that $T\supset[2]$ does not exist}, a lattice $T\supset\bU(2)$
exists if and only if
\roster*
\item
for each prime $p>2$, one has either $\ell(\CN_p)\le r-3$ or $\ell(\CN_p)=r-2$
and, in the latter case,
$\det\CN_p=-\ls|\CN|\bmod(\Z_p\units)^2$;
\item
one has $\CN_2\simeq\Cal{U}_2\oplus\CN_2'$, where $\Cal U_2:=\discr\bU(2)$;
in other words, there is a pair $u,v\in\CN_2$ such that
$u^2=v^2=0\bmod2\Z$ and
$u\cdot v=\frac12\bmod\Z$.
\endroster

\section{Proof of \autoref{th.main}}\label{proof.main}

By abuse of language we call a configuration~$\graph$ of lines/$h$-fragments
\emph{totally real} if $N$, $T$, $\Gf$ in~\eqref{eq.phi} and
$g\in\eqref{eq.real}$ can be \emph{chosen} so that $g|_N=-\id$, \ie, in other words,
there \emph{exists} a polarized $K3$-surface with the given graph of lines and a real
structure with respect to which \emph{all} lines are real.

One of our main observations is the fact that the vast majority of the
graphs found in~\cite{degt:hplanes} are totally real; this is
easily established using \autoref{s.totally.real}, and we merely state the
result. Therefore, we engage into the more involved analysis of the other
real structures as in \autoref{s.real} only when we do not immediately get
$\maxRR(2d)=\maxC(2d)$ and extra work is needed to find the intermediate
bound $\maxR(2d)$.

\subsection{Quartics \rm($h\sp2=4$)}\label{s.h=4}
In the \emph{known} examples ($\ls|\Fn X|\ge44$, but the list is known to be
complete only for $\ls|\Fn X|>48$, see
\cite{degt.Rams:quartics,degt.Rams.anc}),
\roster*
\item
the numbers of $h$-fragments are $20\dt34,36,40\dt42,48,50,60,72$, and
\item
those of totally real $h$-fragments are $22\dt34,36,40,42,48,50$.
\endroster
The maximal \emph{known} number~$50$ is attained at the real surface $Y_{56}$
with $56$ real lines, see \cite[Remark~3.2]{degt:hplanes}. For the other
surfaces with many $h$-fragments, we have:
\roster*
\item
Schur's quartic~$X_{64}$, with the transcendental lattice $T\simeq[8,4,8]$,
has a real structure with $16$ real $h$-fragments
(composed of $28$ real lines, which is the maximum for $X_{64}$,
see \cite[Proposition 9.1]{DIS});
\item
the quartic $X'_{60}$, $T\simeq[4,2,16]$, has a real structure with $12$
real $h$-fragments (composed of $34$ lines, of which only $10$ are real);
\item
the surface $X''_{60}$, $T\simeq[4,1,14]$, has no real structure at all:
that is why, in all statements, we mention a pair $X''_{60},\bar X''_{60}$
of conjugate quartics.
\endroster

\subsection{Sextics \rm($h\sp2=6$)}\label{s.h=6}
Recall that the $h$-fragments in degree~$6$ are
\roster
\item\label{h=6,g=1}
  the $3$-prism, \ie, the (unique) $3$-regular union $K(3)\cup K(3)$, and
\item\label{h=6,g=2}
  the complete bipartite graph $K(3,3)$.
\endroster
Out of the $9235$ $h$-configurations, there are but $357$ that are not
totally real. They have few (at most $30$) $h$-fragments, and
we do not
consider them.

The counts taken by the number of totally real $h$-fragments are shown in
\autoref{tab.h=6}.

\table
\caption{Real $h$-configurations of degree~$6$ \noaux{(see \autoref{conv.tab1})}}\label{tab.h=6}
\vspace{-2\bigskipamount}
\[*
\let\1\tabconfig
\let\u\cross
\minitab{\kern6pt}
\HHH&(r,g,s)&\omit\hss totally real $h$-fragment counts\hss
                       &&\multispan5\minitabskip\hss large counts\hss\\\trule
\iref{h=6,g=1}&(6,3,12) &\clist{0\dt6, 8\dt10, 12, 15, 16, 18, 20, 25, 36}
                       &&16&\.&36&25&36\\
\iref{h=6,g=2}&(6,4,72) &\clist{0\dt24, 26, 27, \u{29}, 30, 36, 40}
                       &&20&36&12&24&40\\
\counts{total counts}   &\clist{0\dt34, 36, 48, 49, 76}
                       &&36&36&48&49&76\\
\counts{configuration}& &&34&\1{\Theta_{36}}&\1{\Psi_{36}'}&\1{\Psi_{38}}&\1{\Psi_{42}}
\endminitab
\]
\endtable

\convention\label{conv.tab1}
Listed in Tables~\ref{tab.h=6} and \ref{tab.h=8}
are the values taken by the number of \emph{totally real} $h$-fragments,
both total and itemized by the type of the graphs.
(The graphs are indexed according to~\cite{degt:hplanes}. For
identification, we present the values of $r:=\rank\Z\HH$, $g:=\girth(\HH)$, and
$s:=\ls|\Aut\HH|$.)
For better comparison, we
also show the complex counterparts of these counts as in~\cite{degt:hplanes}:
the very few values that do not appear in the real
world are struck out.

The last column contains
the itemized counts for a few largest $h$-configurations. In the last row,
for each of these largest $h$-configurations, we show the corresponding (minimal)
configuration of lines, either by its ``name'' in the notation of
\cite{degt:lines,DR:octic.graphs,degt.Rams:octics,degt.Rams:sextics}
or just by the number of
lines.
\endconvention

The three largest counts are realized by triangular, hence special,
configurations $\Psi_{42}$, $\Psi_{38}$, $\Psi'_{36}$. The next one is
$\Theta_{36}$, which is non-special (see
\cite{degt.Rams:sextics}). Together, these observations
prove~\eqref{eq.6nspec}.

Next, recall (see
\cite[Corollary~6.9]{degt.Rams:sextics}) that a Fano graph realized by both
non-special and special smooth sextics must be bipartite of size at most $12+12=24$; in
particular, the $h$-fragments in such a graph are all of type $K(3,3)$.
Found in~\cite{degt:hplanes} are all configurations with $16$ or more
$h$-fragments, and only one of them, \viz. that of Humbert sextics, is
bipartite. This configuration does have $16$ $K(3,3)$-fragments, it is
realized by both non-special and special (a codimension~$1$ stratum) sextics
(see \cite[Remark~7.10]{degt.Rams:sextics}), and, by~\eqref{eq.totally.real} and
\cite[Theorem~4.6]{DDK}, it is totally real. Thus, we have~\eqref{eq.6both}
and the Humbert sextic is the only one realizing these values.

\subsection{Octics \rm($h\sp2=8$)}\label{s.h=8}
Recall that the $h$-fragments in degree~$8$ are
\roster
\item\label{h=8,g=1}
  the (unique) $3$-regular union $K(3)\cup K(3,2)$,
\item\label{h=8,g=2}
  the Wagner graph, and
\item\label{h=8,g=3}
  the $1$-skeleton of a $3$-cube.
\endroster
The former appears on special octics only, whereas the two latter appear only on
triquadrics, see \cite[Proposition~5.1]{degt:hplanes}. In particular, unlike
the case of sextics, we do not need to consider the ``common'' bounds: they
are trivial.

There are but $15$ (out of $860$) $h$-configurations
that are not totally real; all except two have few (at most $24$) $h$-fragments. The
values taken by the $h$-fragment counts are given in \autoref{tab.h=8},
where $\HH_1$ and $\HH_2,\HH_3$ are treated separately.

\table
\caption{Real $h$-configurations of degree~$8$ \noaux{(see \autoref{conv.tab1})}}\label{tab.h=8}
\vspace{-2\bigskipamount}
\[*
\let\1\tabconfig
\let\u\cross
\minitab{\kern6pt}
\HHH&(r,g,s)&\omit\hss totally real $h$-fragment counts\hss
                       &&\multispan6\minitabskip\hss large counts\hss\\\trule
\iref{h=8,g=1}&(8,3,12) &\clist{0\dt10, 12, 14, 16, 18, 20, 36, 42, \u{72}}
                       &&42\\
\iref{h=8,g=2}&(8,4,16) &\clist{0\dt16, 18\dt21, 24, 32, 48}
                       &&&24&16&48&32&\.\\
\iref{h=8,g=3}&(8,4,48) &\clist{0\dt13, 16, 20, 21, 32, 80}
                       &&&12&20& 8&32&80\\
\counts{total $\iref{h=8,g=2}+\iref{h=8,g=3}$}
                        &\clist{0\dt18, 20, 21, 23, \u{24}, 26, 36, 56, 64, 80}
                       &&&36&36&56&64&80\\
\counts{configuration}&&&\1{\Psi_{29}}&28&28&\1{\Theta_{32}'}&\1{\Theta_{32}''}&\1{\Theta_{32}^\mathrm{K}}
\endminitab
\]
\endtable

Missing from \autoref{tab.h=8} (as compared to \cite[Table~2]{degt:hplanes})
are two large $h$-configurations
$\Psi_{33}$ with $72\times\HH_1$ and $\Theta'''_{32}$ with
$(48\times\HH_2)\cup(16\times\HH_3)$. Both have $T\simeq2\bU(3)$,
and
\roster*
\item
$\Psi_{33}$ has a real form with $54$ real $h$-fragments;
\item
the real forms of $\Theta'''_{32}$ have at most $8$ real $h$-fragments.
\endroster
In particular, this gives us a proof of~\eqref{eq.8spec}.

\subsection{The polarization $h\sp2=10$}\label{s.h=10}
There are $13$ $h$-configurations that are not totally real. The two large
ones have $T\simeq[4,2,10]$:
\roster*
\item
$(6\times\HH_2)\cup(6\times\HH_3)$:
there are no real forms with real $h$-fragments;
\item
$14\times\HH_6\simeq\Phi'''_{30}$ in \cite{degt:lines}:
there are real forms with $2$, $4$, or~$6$ real $h$-fragments.
\endroster
The other eleven have at most nine $h$-fragments; we do not consider them.

\subsection{The polarization $h\sp2=14$}\label{s.h=14}
The three $h$-configurations that are not totally real are
\roster*
\item
$3\times\HH_4$, with $T\simeq\bA_2(3)$
and real forms with a single real $h$-fragment,
\item
$4\times\HH_6$, with $T\simeq2\bU(3)$
and real forms with two real $h$-fragments, and
\item
$2\times\HH_8$, with $T\simeq\bA_2(3)$
and real forms with two real $h$-fragments.
\endroster

\subsection{The polarization $h\sp2=16$}\label{s.h=16}
The only $h$-configuration that is not totally real is
$\Delta'_{32}$ (see~\cite{degt:lines}) with $24$ copies of $\HH_8$: it has
$T\simeq\bA_2(2)$ and admits a real form with four or six real $h$-fragments.
In the latter case, the surface does not have a single real line: the real
structure interchanges the two Kummer $16$-tuples.
It is this real form that realizes the upper bound $\maxR(16)=\underline6$.

\subsection{Other polarizations}\label{s.other}
If $h^2\in\{12,18,20,22,24,28\}$, all $h$-configurations are totally real.
This completes the proof.
\qed

\section{Hyperelliptic models}\label{S.hyper}

For hyperelliptic surfaces, \ie, two-to-one maps $X\to Y\subset\Cp{d+1}$ onto
a rational surface~$Y$ with smooth ramification locus, an analogue of the table in
\autoref{th.main} is as follows:
\[*
\let\u\underline
\minitab\quad
h^2=2d          &     2&  4& 6& 8\\\trule
\maxC           &    72&144&36&56\\
\maxR           &\u{66}&144&36&56\\
\maxRR          &\u{66}&144&36&56\\
\endminitab
\]

For sextics ($h^2=6$) and octics ($h^2=8$),
all $h$-configurations are totally real; hence, as in the
complex case, the values taken are
\[*
\numRR(X)\in
\begin{cases}
  \{0, 1, 3, 6, 10, 15, 21, 28, 36\}, &\text{if $2d=6$}, \\
  \{0, 1, 4, 10, 20, 56\}, &\text{if $2d=8$}.
\end{cases}
\]
As in the other cases, we do not engage into the meticulous analysis of other
real structures on these surfaces.

\subsection{Double planes \rm($h\sp2=2$)}\label{s.P2}
All $2$-polarized $K3$-surfaces
(which, for the sake of completeness, also appear in \autoref{th.main}) are
double planes $X\to\Cp2$ ramified over smooth sextic curves $C\subset\Cp2$. A
real structure $\iota\:\Cp2\to\Cp2$ preserving~$C$ lifts to two
 real structures $\iota',\iota''\:X\to X$ that differ by the deck
translation of the double covering. Therefore, the real $h$-fragments on~$X$
are the pull-backs of the real tritangents to~$C$; the two lines of such a
pull-back are both real with respect to one of $\iota',\iota''$, and
they are complex conjugate with respect to the other. However, since any two
tritangents intersect in $\Cp2$, we conclude that there is a common lift,
say, $\iota'$, that makes all real pull-backs totally real simultaneously.
Thus, depending on the choice of $\iota',\iota''$,
\[*
\numRR(X)=0\quad\text{or}\quad\numRR(X)=\numR(X).
\]
Using the description of real tritangents found in~\cite{degt:sextics}), we
arrive at
\[*
\maxRR(2)=\maxR(2)=\underline{66}\quad
 \text{(\vs. $\maxC(2)=72$, see \cite{degt:hplanes})},
\]
but the only known values taken by $\maxRR(2)$ are $63$ and~$66$.

\subsection{Double quadrics \rm($h\sp2=4$)}\label{s.2x2}
The configuration
maximizing the
number of $h$-fragments is totally real; hence,
\[*
\maxRR(4\hyp)=\maxR(4\hyp)=\maxC(4\hyp)=144.
\]

For some reason, lines on hyperelliptic quartics (equivalently, generatrices
that are bitangent to a bidegree~$(4,4)$ curve on a quadric $\Cp1\times\Cp1$) have not
be studied in detail in~\cite{degt:lines}. We bridge this gap, starting from
a certain number $1\le n\le12$ of generatrices in one family and using the
line-by-line algorithm of \cite[\S\,2.7]{degt.Rams:sextics} to add a number
$0\le m\le n$ of bitangents in the other family.
The two extreme cases are
\roster*
\item
if $m=0$, then $n=0\dt7,8,8,9,9,10,10,11,12$ ($16$ strata, all totally real);
\item
if $n=12$, then $m=0,4,8,12$ (four strata, all totally real).
\endroster
More generally, the pair $m\le n$ takes the following values:
\[*
\minitab\quad
n=&0\dt 9&      10&       11&      12\\
m=&0\dt n&0\dt8,10&0,2,4,6,8&0,4,8,12
\endminitab
\]

There are $470$ $h$-configurations, of which only $16$ are not totally real.
The values taken by the $h$-fragment count $mn$, both complex and totally real, are
\begin{multline*}
0\dt10, 12, 14\dt16, 18, 20\dt22, 24, 25, 27, 28, 30, 32, 35, 36, 40, 42, 44,\\
45, 48\dt50, 54, 56, 60, 63, 64, 66, 70, 72, 80, 81, 88, 96, 100, 144.
\end{multline*}
(Recall that the $h$-fragments on a double quadric $X\to\Cp1\times\Cp1$
are the pull-backs of the pairs
of intersecting bitangents to the ramification locus.)

{
\let\.\DOTaccent
\def\cprime{$'$}
\bibliographystyle{amsplain}
\bibliography{degt}

\providecommand{\bysame}{\leavevmode\hbox to3em{\hrulefill}\thinspace}
\providecommand{\MR}{\relax\ifhmode\unskip\space\fi MR }
\providecommand{\MRhref}[2]{%
  \href{http://www.ams.org/mathscinet-getitem?mr=#1}{#2}
}
\providecommand{\href}[2]{#2}
\begin{thebibliography}{10}

\bibitem{degt:lines}
Alex Degtyarev, \emph{Lines on {S}mooth {P}olarized {K}3-{S}urfaces}, Discrete
  Comput. Geom. \textbf{62} (2019), no.~3, 601--648. \MR{3996938}

\bibitem{degt:sextics}
\bysame, \emph{Tritangents to smooth sextic curves}, Ann. Inst. Fourier
  (Grenoble) \textbf{72} (2022), no.~6, 2299--2338. \MR{4500357}

\bibitem{degt:hplanes}
\bysame, \emph{Real split hyperplane sections on smooth polarized
  {$K3$}-surfaces}, To appear, \arXiv{2509.24349}, 2025.

\bibitem{DDK}
Alex Degtyarev, Igor Dolgachev, and Shigeyuki Kond{\=o}, \emph{{K3} surfaces of
  degree six arising from desmic tetrahedra}, To appear in Rendiconti Lincei.
  Matematica e Applicazioni, 2025.

\bibitem{DI:empty}
Alex Degtyarev and Ilia Itenberg, \emph{Real plane sextics without real
  points}, J. Algebraic Geom. \textbf{34} (2025), no.~3, 543--577. \MR{4895154}

\bibitem{DIS}
Alex Degtyarev, Ilia Itenberg, and Ali~Sinan Sert\"oz, \emph{Lines on quartic
  surfaces}, Math. Ann. \textbf{368} (2017), no.~1-2, 753--809. \MR{3651588}

\bibitem{DR:octic.graphs}
Alex Degtyarev and S{\l}awomir Rams, \emph{Extended {F}ano graphs of octic
  {$K3$}-surfaces}, Electronic,
  \url{http://www.fen.bilkent.edu.tr/~degt/papers/octics_Fano.zip}, March 2021.

\bibitem{degt.Rams.anc}
\bysame, \emph{Ancillary files for the paper: {L}ines on {$K3$}-quartics via
  triangular sets}, ancillary files for the preprint \arXiv{2301.04127}, 2022.

\bibitem{degt.Rams:octics}
\bysame, \emph{Counting lines with {V}inberg's algorithm}, To appear in Rev.
  Mat. Iberoam., 2025.

\bibitem{degt.Rams:quartics}
\bysame, \emph{Lines on {$K3$}-quartics via triangular sets}, Discrete Comput.
  Geom. \textbf{73} (2025), no.~3, 785--814. \MR{4880184}

\bibitem{degt.Rams:sextics}
\bysame, \emph{Lines on sextic {$K3$}-surfaces with simple singularities}, To
  appear, 2025.

\bibitem{GAP4.13}
The GAP~Group, \emph{{GAP -- Groups, Algorithms, and Programming, Version
  4.13.0}}, 2024.

\bibitem{Itenberg}
I.~V. Itenberg, \emph{Rigid isotopy classification of curves of degree {$6$}
  with one nondegenerate double point}, Topology of manifolds and varieties,
  Adv. Soviet Math., vol.~18, Amer. Math. Soc., Providence, RI, 1994,
  pp.~193--208. \MR{1296896}

\bibitem{Josi}
Johannes Josi, \emph{Nodal rational sextics in the real projective plane},
  Ph.D. thesis, 2018.

\bibitem{Nikulin:forms}
V.~V. Nikulin, \emph{Integer symmetric bilinear forms and some of their
  geometric applications}, Izv. Akad. Nauk SSSR Ser. Mat. \textbf{43} (1979),
  no.~1, 111--177, 238, English translation: Math USSR-Izv. 14 (1979), no. 1,
  103--167 (1980). \MR{525944 (80j:10031)}

\bibitem{Nikulin:sublattice}
\bysame, \emph{Involutions of integer quadratic forms and their applications to
  real algebraic geometry}, Izv. Akad. Nauk SSSR Ser. Mat. \textbf{47} (1983),
  no.~1, 109--188. \MR{688920}

\end{thebibliography}
}

\end{document}